\documentclass[12pt, a4paper, reqno]{amsart}

\usepackage[english]{babel}

\usepackage[T1]{fontenc}
\usepackage{lmodern}

\usepackage{amsmath}
\usepackage{amssymb}
\usepackage{amsthm}

\usepackage{url}

\usepackage{hyperref} 

\usepackage{enumitem}
\usepackage[normalem]{ulem}

\renewcommand{\epsilon}{\varepsilon}
\renewcommand{\theta}[0]{\vartheta}
\renewcommand{\phi}[0]{\varphi}

\newcommand{\gammabar}[0]{\bar{\gamma}}

\newcommand{\ppar}{\ \par}

\newcommand{\Span}[1]{\left\langle\, #1 \,\right\rangle}

\newcommand{\Set}[1]{\left\{ #1 \right\}}

\newcommand{\abuse}[1]{[\![ #1 ]\!]}

\DeclareMathOperator{\Aut}{Aut}

\DeclareMathOperator{\Ann}{Ann}

\DeclareMathOperator{\Hol}{Hol}

\newtheorem{dummy}{Dummy}
\numberwithin{dummy}{section}
\numberwithin{figure}{section}

\newtheorem{theorem}[dummy]{Theorem}

\newtheorem*{theoremstar}{Theorem}

\theoremstyle{definition}

\theoremstyle{remark}

\newtheorem{remark}[dummy]{Remark}

\makeatletter
\def\imod#1{\allowbreak\mkern10mu({\operator@font mod}\,\,#1)}
\makeatother

\numberwithin{equation}{section}

\usepackage{amsthm, mathtools, amssymb, amsfonts, stmaryrd, latexsym, mathrsfs, todonotes, eufrak, color}
\usepackage[all]{xy}
\usepackage{tikz-cd}
\usetikzlibrary{arrows}

\begin{document}

\date{23 June 2023, 18:26 CEST --- Version 0.16%
}

\title[Isoclinism of skew braces]
      {Some properties of skew braces
      \\
      that are invariant under isoclinism}
      
\author{A. Caranti}

\address[A.~Caranti]%
 {Dipartimento di Matematica\\
  Universit\`a degli Studi di Trento\\
  via Sommarive 14\\
  I-38123 Trento\\
  Italy} 

\email{andrea.caranti@unitn.it} 

\urladdr{https://caranti.maths.unitn.it/}



\begin{abstract}
  Letourmy and Vendramin have recently introduced
  a concept of isoclinism for skew braces.
  
  We show that for a skew brace the properties of being bi-skew,
  $\lambda$-homomorphic, and inner are invariant under isoclinism.
\end{abstract}

\thanks{The author is a member of INdAM---GNSAGA. The  author
  gratefully acknowledges support from the Department of Mathematics of
  the University of Trento.}

\maketitle

\thispagestyle{empty}

\section{Introduction}

In 1940 Philip Hall introduced a relation weaker than
isomorphism between groups, called isoclinism~\cite{Hall}. Two groups
are said to 
be isoclinic if there is a pair of isomorphisms between their central
quotients, resp.\ commutator subgroups, which are compatible with the
commutator map. (In Sections~\ref{sec:ODNT}~and
\ref{sec:iso} we will give full definitions for the concepts
dicussed in this Introduction.)

In 2022 Thomas Letourmy and Leandro Vendramin have introduced a
concept of isoclinism for skew braces~\cite{LV}. At the ``Advances 
in Group Theory and Generalizations'' Conference held in Lecce, Italy,
in June 2023, Vendramin asked whether the concept of being a bi-skew
brace is  invariant under isoclinism --- bi-skew braces have been
introduced by Lindsay Childs in~\cite{bi-Childs}.

At the same Conference, Lorenzo  Stefanello suggested to consider also
the behaviour  under isoclinism of $\lambda$-homomorphic  skew braces,
as defined  by Valeriy Bardakov,  Mikhail Neshchadim, and  Manoj Yadav
in~\cite{lambda}, and of inner skew braces.

The goal of this short note is to show that all these concepts are
indeed invariant under isoclinism:
\begin{theoremstar}
  For skew braces, the properties of being
  \begin{enumerate}
  \item a bi-skew brace,
  \item a $\lambda$-homomorphic skew brace, and
  \item an inner skew brace
  \end{enumerate}
  are invariant under isoclinism.
\end{theoremstar}

We set up notation in Section~\ref{sec:ODNT}. In Section~\ref{sec:iso}
we discuss isoclinism. Section~\ref{sec:bi} deals with bi-skew braces,
while  Section~\ref{sec:other} deals  with $\lambda$-homomorphic  skew
braces and inner skew braces.

We are grateful to Leandro
Vendramin and Lorenzo Stefanello for suggesting these problems.

\section{Old Dogs and New Tricks}
\label{sec:ODNT}

In the literature a (left) skew brace is usually defined as a
triple $(B, +, 
\circ)$, where $(B, +)$ and $(B, \circ)$ are groups, neither of which
is assumed to be abelian, and the two
operations are related by
\begin{align}
  \label{eq:lsb}
  - x + (x \circ (y + z)) = - x + (x \circ y) - x + (x \circ z),
\end{align}
for $x, y, z \in B$. Such a skew brace can be defined in terms of the
maps, for $x \in B$,
\begin{equation*}
  \label{eq:lambda}
  \begin{aligned}
    \lambda_{x} :\ &B \to B\\
    &y \mapsto -x + (x \circ y)
  \end{aligned}
\end{equation*}
which~\eqref{eq:lsb} shows to be  endomorphisms of $(B, +)$.

However,  this author  has long  been used  to writing  endomorphisms as
exponents (which does not fit well  with an additive notation), and to
using actions on  the right (whereas the  setting above implies actions
on the left).
Therefore in  this paper we will  abuse the patience of  the reader by
employing the following notation, which stems from a related context
dealt with in~\cite{CDVS}.

A \emph{(right) skew brace} will be a triple $(B, \cdot, \circ)$, where 
$(B, \cdot)$ and $(B, \circ)$ are groups, and the two
operations are related by
\begin{align*}
  \label{eq:rsb}
  ((y \cdot z) \circ x) \cdot x^{-1}
  =
  (y \circ x) \cdot x^{-1}
  \cdot
  (z \circ x) \cdot x^{-1}
\end{align*}
for $x, y, z \in B$, where
$x^{-1}$ denotes the inverse in $(B, \cdot)$. In this notation, this
states that for each $x \in B$ the function
\begin{align*}
  \gamma(x) :\ &B \to B
  \\
               &y \mapsto   (y \circ x) \cdot x^{-1} 
\end{align*}
is an endomorphism of $(B, \cdot)$. In fact, each $\gamma(x)$ is an
automorphism of $(B, \cdot)$, and the
\emph{gamma function} $\gamma : (B, \circ) \to
\Aut(B, \cdot)$ is a morphism of 
groups, that is, it satisfies the functional equation
\begin{equation}
  \label{eq:gfe}
  \gamma(y^{\gamma(x)} \cdot x) = \gamma(y) \gamma(x),
\end{equation}
where the operation  on the right is composition  of maps. Conversely,
if $(B, \cdot)$ is  a group, and $\gamma : B \to  \Aut(B, \cdot)$ is a
map that satisfies~\eqref{eq:gfe}, then  defining ``$\circ$''  as $y
\circ x =  y^{\gamma(x)} \cdot x$, we obtain a  skew brace $(B, \cdot,
\circ)$.

In keeping with the more common nomenclature, we will refer to $(B,
\cdot)$ as the \emph{additive} group and to $(B, \circ)$ as the
\emph{multiplicative} group of the skew brace.

\section{Isoclinism}
\label{sec:iso}

For a group $B = (B, \cdot)$ and $x, y \in B$, we write the product $x
\cdot y$ as juxtaposition $x y$, unless clarity suggests otherwise. We
also write $x^{y} =
y^{-1} x y$ for the conjugate of $x$ by $y$, $[x, y] = x^{-1} y^{-1} x y = x^{-1}
x^{y}$ for the 
commutator of $x$ and $y$,
\begin{equation*}
  [B, B]
  =
  \Span{ [x, y] : x, y \in G }
\end{equation*}
for the commutator subgroup of $G$, and 
\begin{equation*}
Z(B, \cdot) = \Set{ z \in B : [x, z] = 1
  \text{ for all $x \in B$} }
\end{equation*}
for its centre. Note that $x^{y} = x
\cdot [x, y]$. We record 
the standard commutator identities
\begin{equation}
  \label{eq:standard}
  \begin{aligned}
    &
    [x, y z]
    =
    [x, z] \cdot [x, y]^{z}
    =
    [x, z] \cdot [x, y] \cdot [[x, y], z],
    \\&
    [x y, z]
    =
    [x, z]^{y} \cdot [x, y]
    =
    [x, z] \cdot [[x, z], y] \cdot [x, y].
  \end{aligned}
\end{equation}

We  now review  in  our notation  the concept  of  isoclinism of  skew
braces, as introduced in~\cite{LV}.

Let $(B, \cdot, \circ)$ be a skew brace. It is useful to compute in
the holomorph $\Hol(B, \cdot) = \Aut(B, \cdot)
\ltimes (B, \cdot)$ of the group $(B, \cdot)$. A natural replacement
for the centre is the annihilator 
\begin{equation*}
  \Ann(B)
  =
  Z(B, \cdot)
  \cap
  \ker(\gamma)
  \cap
  C_{B}(\gamma(B)),
\end{equation*}
where $C_{B}(\gamma(B)) = \Set{ x \in B : x^{\gamma(y)} = x \text{ for
    all $y \in B$} }$ denotes the centralizer in $B$ of $\gamma(B)$
within $\Hol(B, \cdot)$. The annihilator $\Ann(B)$ is an ideal of the
skew brace $B$, so that one may consider the quotient skew brace
$(B / \Ann(B), \cdot, \circ)$. Besides the \emph{additive commutator} $[x,
  y] = x^{-1} y^{-1} x y$ of $x, y \in B$, consider the expression
\begin{equation}
  \label{eq:star}
  \begin{aligned}
    x * y
    =
    x^{-1} \cdot (x \circ y) \cdot y^{-1}.
  \end{aligned}
\end{equation}
We call this the \emph{$*$-commutator} because the calculation
\begin{equation}
  \label{eq:starcomm}
  x * y
  =
  x^{-1} \cdot x^{\gamma(y)} \cdot y \cdot y^{-1}
  =
  x^{-1} \cdot x^{\gamma(y)}
  =
  [x, \gamma(y)],
\end{equation}
shows that this is indeed  a bona fide commutator in $\Hol(B,
\cdot)$. Both expressions~\eqref{eq:star}~and \eqref{eq:starcomm} will
be useful in the sequel.

It is proved
in~\cite[Proposition 2.2]{LV} that the additive subgroup
\begin{equation*}
  B'
  =
  [B, B] \cdot [B, \gamma(B)]
\end{equation*}
of $B$ is an ideal of $B$.

\eqref{eq:standard} shows that the map
\begin{equation}
  \label{eq:commcomm}
  \begin{aligned}
    \abuse{\cdot, \cdot}
    :\ &B / \Ann(B) \times B / \Ann(B) \to [B, B]\\
    &(x \Ann(B), y \Ann(B)) \mapsto [x, y]
  \end{aligned}
\end{equation}
is well defined,  as $\Ann(B) \subseteq Z(B, \cdot)$; in  fact,
regarding
$x \Ann(B)$ and $y \Ann(B)$ as subsets of $B$, we have
\begin{equation*}
  [x \Ann(B), y \Ann(B)]
  =
  \Set{ [x c, y d] : c, d \in \Ann(B) }
  =
  \Set{ [x, y] }.
\end{equation*}
We will thus write the map~\eqref{eq:commcomm} as
\begin{equation}
  \label{eq:abuse1}
  \abuse{x \Ann(B), y \Ann(B)} = [x, y].
\end{equation}

Let
\begin{equation*}
  \gammabar : (B / \Ann(B), \circ) \to \Aut(B / \Ann(B), \cdot)
\end{equation*}
be the gamma function of $B / \Ann(B)$.
We have
\begin{align*}
  (x \Ann(B))^{\gammabar(y \Ann(B))}
  &=
  (x \Ann(B) \circ y \Ann(B)) \cdot (y \Ann(B))^{-1}
  \\&=
  ((x \circ y) \cdot y^{-1}) \Ann(B)
  \\&=
  x^{\gamma(y)} \Ann(B).
\end{align*}
In particular,
\begin{equation}
  \label{eq:ataleoftwgammas}
  [x \Ann(B), \gammabar(y \Ann(B))]
  =
  [x, \gamma(y)] \Ann(B).
\end{equation}
Moreover, $\gamma$ is constant on the additive
cosets $x \Ann(B)$ of $\Ann(B)$: if $y \in B$ and $d \in \Ann(B)$, we
have
\begin{equation*}
  \gamma(y d)
  =
  \gamma(y^{\gamma(d)} d)
  =
  \gamma(y) \gamma(d)
  =
  \gamma(y),
\end{equation*}
as $\Ann(B) \subseteq \ker(\gamma)$.
To show that the $*$-commutator also yields a well
defined map
\begin{equation}
  \label{eq:star-is-well-defined}
  \begin{aligned}
    * :\ &B / \Ann(B) \times B / \Ann(B) \to B * B = [B, \gamma(B)]\\
    &(x \Ann(B), y \Ann(B)) \mapsto x * y = [x, \gamma(y)]
  \end{aligned}
\end{equation}
it remains to prove that for $x, y \in B$ and
$c \in \Ann(B)$
one has $(x c) * y = x * y$. This follows from
\begin{equation*}
  (x c) * y
  =
  [x c , \gamma(y)]
  =
  [x, \gamma(y)]^{c} \cdot [c, \gamma(y)]
  =
  [x, \gamma(y)]
  =
  x * y,
\end{equation*}
as $\Ann(B) \subseteq Z(B, \cdot) \cap C_{B}(\gamma(B))$. As above, we will
write the map~\eqref{eq:star-is-well-defined} as
\begin{equation}
  \label{eq:abuse2}
  \abuse{x \Ann(B), \gammabar(y \Ann(B))}
  =
  [x, \gamma(y)],
\end{equation}
in a notation that we have just seen to be well defined.

Let $(B_{1}, \cdot, \circ), (B_{2}, \cdot, \circ)$ be two skew
braces. We will write $\gamma$ for the gamma functions of both, and
similarly for the ($*$-)commutator and for $\abuse{\cdot,\cdot}$.
An \emph{isoclinism} between the two is a pair of skew brace isomorphisms 
\begin{equation*}
  \xi : B_{1} / \Ann(B_{1}) \to B_{2} / \Ann(B_{2}),
  \quad
  \theta : B_{1}' \to B_{2}'
\end{equation*}
such that the following diagram commutes.
\begin{equation}
  \label{eq:isodiagram}
  \begin{tikzcd}[column sep = 1.5cm, row sep = 1.5cm]
    B_{1}'
    \arrow{d}[swap]{\theta}
    &
    \arrow{l}[swap]{*}
    B_{1} / \Ann(B_{1}) \times B_{1} / \Ann(B_{1})
    \arrow{d}{\xi \times \xi}
    \arrow{r}{{\abuse{\cdot,\cdot}}} & {B_{1}'}
    \arrow{d}{\theta}\\
    B_{2}'
    &
    \arrow{l}{*}
    B_{2} / \Ann(B_{2}) \times B_{2} / \Ann(B_{2})
    \arrow[swap]{r}{{\abuse{\cdot,\cdot}}} & {B_{2}'}
  \end{tikzcd}
\end{equation}
In other words, writing $A_{i} = \Ann(B_{i})$ for $i = 1, 2$, we have
for $x, y \in B_{1}$ 
\begin{equation}
  \label{eq:isoclinism}
  \begin{aligned}
    &
    \abuse{x A_{1}, y A_{1}}^{\theta}
    =
    \abuse{(x A_{1})^{\xi}, (y A_{1})^{\xi}},
    \\
    &
      \abuse{x A_{1}, \gammabar(y A_{1})}^{\theta}
      =
      \abuse{(x A_{1})^{\xi}, \gammabar((y A_{1})^{\xi})},
  \end{aligned}
\end{equation}

\begin{remark}\ppar
  \begin{enumerate}
  \item
    Taken by itself, the first identity of~\eqref{eq:isoclinism}, that
    is,    the    fact    that    the   rightmost    part    of    the
    diagram~\eqref{eq:isodiagram} commutes, defines a slight variation
    on the  concept of isoclinism  of groups, as introduced  by Philip
    Hall in~\cite{Hall}. Hall's version has  $Z(B, \cdot)$ in place of
    $\Ann(B)$.
  \item 
    In a skew brace $(B, \cdot, \circ)$ there is another natural notion of
    a \emph{multiplicative commutator} of $x, y \in B$, namely
    \begin{equation*}
      [x, y]_{\circ}
      =
      x^{\ominus 1} \circ y^{\ominus 1} \circ x \circ y,
    \end{equation*}
    where $x^{\ominus 1}$ denotes the inverse in $(B, \circ)$.
    It is
    proved in~\cite[Proposition 2.15]{LV} that if the two skew braces
    $(B_{1}, \cdot, \circ), (B_{2}, \cdot, \circ)$ are isoclinic, then,
    in the above notation,
    $[B, B]_{\circ} \subseteq B'$, 
    the multiplicative
    commutator is well defined on $B_{i} / A_{i}$, and
    for all $x, y \in B_{1}$ one has
    \begin{equation*}
      \abuse{x  A_{1},   y  A_{1}}_{\circ}^{\theta}
      =
      \abuse{(x   A_{1})^{\xi},  (y
        A_{1})^{\xi}}_{\circ}.
    \end{equation*}
  \end{enumerate}
\end{remark}

\section{Bi-skew braces}
\label{sec:bi}

A \emph{bi-skew brace}, as defined by Childs in~\cite{bi-Childs}, is a skew
brace $(B, \cdot, \circ)$ such that $(B, \circ, \cdot)$ is also a skew brace.
Bi-skew braces have been studied by the author in~\cite{mybi}, and by
Stefanello and Trappeniers in~\cite{ST}. We will use the following
characterisation
\begin{theorem}[\protect{\cite[Theorem 3.1]{mybi}}]
  \ppar
  For a skew brace $(B, \cdot, \circ)$, the following are equivalent:
  \begin{enumerate}
  \item $B$ is a bi-skew-brace, 
  \item the gamma function is an anti-homomorphism $(B, \cdot) \to
    \Aut(B, \cdot)$, that is
    $
      \gamma(y z) = \gamma(z) \gamma(y)
    $
    for all $y, z \in B$, and
  \item
    for all $x, y, z \in B$ one has
    $
      [x, \gamma(y z)]
      =
      [x, \gamma(z) \gamma(y)].
    $
  \end{enumerate}
\end{theorem}

Keeping the notation of the previous paragraph, let $(B_{1}, \cdot,
\circ), (B_{2}, \cdot, \circ)$ be two isoclinic skew braces. Assume
$(B_{1}, \cdot, \circ)$ to be bi-skew. We will prove that $(B_{2}, \cdot,
\circ)$ is also bi-skew.

Let $u, v, w \in B_{2}$, and $x, y, z \in B_{1}$ such that
$u A_{2} = (x A_{1})^{\xi}, v A_{2} = (y A_{1})^{\xi}, w A_{2} =
(z A_{1})^{\xi}$. We have, according to~\eqref{eq:standard}~and
\eqref{eq:ataleoftwgammas} 
\begin{equation}
  \begin{aligned}
    \label{eq:B1-is-bi}
          \abuse{x A_{1}, \gammabar(y z A_{1})}
          &=
          [x, \gamma(y z)]
          =
            [x, \gamma(z) \gamma(y)]
            \\&=
              [x, \gamma(y)]
              \cdot
                  [x, \gamma(z)]
                  \cdot
                  [[x, \gamma(z)], \gamma(y)]
            \\&=
              \abuse{x A_{1}, \gammabar(y A_{1})}
              \cdot
                  \abuse{x A_{1}, \gammabar(z A_{1})}
                  \cdot\\&\ \qquad
                      \cdot \abuse{ [ x , \gamma(z)] A_{1},
                        \gammabar(y A_{1})}. 
  \end{aligned}
\end{equation}

Applying $\theta$ to the left hand side of~\eqref{eq:B1-is-bi}, we get
\begin{align}
  \label{eq:LHS}
  \abuse{x A_{1}, \gammabar(y z A_{1})}^{\theta}
  =
  \abuse{u A_{2}, \gammabar(v w A_{2})}
  =
  [u, \gamma(v w)].
\end{align}

Applying $\theta$ to the right hand side of~\eqref{eq:B1-is-bi}, we get
\begin{multline}
  \label{eq:RHS}
  (\abuse{x A_{1}, \gammabar(y A_{1})}
  \cdot
      \abuse{x A_{1}, \gammabar(z A_{1})}
      \cdot 
  \abuse{[x, \gamma(z)] A_{1}, \gammabar(y A_{1})})^{\theta}
  =\\=
  \abuse{u A_{2}, \gammabar(v A_{2})}
  \cdot
      \abuse{u A_{2}, \gammabar(w A_{2})}
      \cdot
      \abuse{ ([ x , \gamma(z)] A_{1})^{\xi}, \gammabar( (y A_{1})^{\xi} )}.
\end{multline}
Since $\xi : B_{1}/A_{1}\to B_{1}/A_{1}$ is an
isomorphism for both operations ``$\cdot$'' and ``$\circ$'' we have,
according to~\eqref{eq:star}~and \eqref{eq:ataleoftwgammas},
\begin{align*}
  ([ x , \gamma(z)] A_{1})^{\xi}
  &=
  ((x * z) A_{1})^{\xi}
  =
  ((x^{-1} \cdot (x \circ z) \cdot z^{-1}) A_{1})^{\xi}
  \\&=
  ((x A_{1})^{-1} \cdot ( x A_{1} \circ z A_{1}) \cdot (z A_{1})^{-1})^{\xi}
  \\&=
  (((x A_{1})^{\xi})^{-1} \cdot
  ( (x A_{1})^{\xi} \circ (z A_{1})^{\xi})
  \cdot ((z A_{1})^{\xi})^{-1})
  \\&=
  (x A_{1})^{\xi} * (z A_{1})^{\xi}
  =
    [(x A_{1})^{\xi}, \gammabar((z A_{1})^{\xi})]
    \\&=
    [u A_{2}, \gammabar( w A_{2})]
    =
    [u , \gamma( w )] A_{2}.
\end{align*}
Thus the last term of the right hand side of~\eqref{eq:RHS} equals
\begin{equation*}
  \abuse{ [u , \gamma( w )] A_{2}, \gammabar(v A_{2})}
  =
  [[u, \gamma(w)], v],
\end{equation*}
so that~\eqref{eq:RHS} equals
\begin{equation}
  \label{eq:RHS2}
        [u , \gamma(v )]
        \cdot
            [u , \gamma(w )]
            \cdot
                [[u , \gamma(w )],
                  \gamma(v )]
                =
                [u, \gamma(w) \gamma(v)].
\end{equation}
Comparing the right hand sides of~\eqref{eq:LHS}~and \eqref{eq:RHS2},
we see that 
$B_{2}$ is a bi-skew brace.

\section{$\lambda$-homomorphic and inner skew braces}
\label{sec:other}

According to the definition in~\cite{lambda}, a skew brace $(B, \cdot,
\circ)$ is said to be \emph{$\lambda$-homomorphic} if the gamma function is a
homomorphism $\gamma : (B, \cdot) \to \Aut(B, \cdot)$, that is,
\begin{equation*}
  \gamma(x y)
  =
  \gamma(x) \gamma(y)
\end{equation*}
for $x, y \in B$. This property had been characterised by Elena Campedel,
Ilaria Del Corso, and the author in~\cite[Lemma 2.13]{p2q}.

The proof that if two skew braces $B_{1}$ and $B_{2}$ are isoclinic,
then $B_{1}$ is $\lambda$-homomorphic if and only if $B_{2}$ is
$\lambda$-homomorphic is an 
obvious variation on the proof of the preceding section.

A skew brace $(B, \cdot,
\circ)$ is said to be \emph{inner} if the gamma function takes values in the
group of inner automorphisms of $(B, \cdot)$. In other words, for $y
\in B$, there is $t \in B$ such that for all $x \in G$ one has
$x^{\gamma(y)} = x^{t}$, or 
equivalently $[x, \gamma(y)] = x^{-1} x^{\gamma(y)} = x^{-1} x^{t} =
[x, t]$. Special classes of inner skew braces have been dealt with
in~\cite{Koch, CS1, CS3}. 

Keeping the notation of Sections~\ref{sec:iso}~and \ref{sec:bi},
assume $B_{1}$ to be inner. Let 
$u, v \in B_{2}$ and $x, y \in B_{1}$ such that 
$u A_{2} = (x A_{1})^{\xi}, v A_{2} = (y A_{1})^{\xi}$.

We have
\begin{align*}
  [u, \gamma(v)]
  &=
  \abuse{u A_{2}, \gammabar(v A_{2})}
  =
  \abuse{(x A_{1})^{\xi}, \gammabar((y A_{1})^{\xi})}
  \\&=
  \abuse{x A_{1}, \gammabar(y A_{1}))}^{\theta}
  =
  [x, \gamma(y)]^{\theta}.
\end{align*}
Since $B_{1}$ is inner, there is $t \in B_{1}$, depending only on $y$,
such that $[x, \gamma(y)] = [x, t]$, so that
\begin{align*}
  [x, \gamma(y)]^{\theta}
  &=
  [x, t]^{\theta}
  =
  \abuse{ x A_{1}, t A_{1}}^{\theta}
  =
  \abuse{ (x A_{1})^{\xi}, (t A_{1})^{\xi} }
  \\&=
  \abuse{ u A_{2}, w A_{2} }
  =
  [u, w],
\end{align*}
where $w \in (t A_{1})^{\xi}$ depends only on $t$, and thus on $y$,
and ultimately on $v$. It follows that $B_{2}$ is also an
inner skew brace.


\providecommand{\bysame}{\leavevmode\hbox to3em{\hrulefill}\thinspace}
\providecommand{\MR}{\relax\ifhmode\unskip\space\fi MR }
\providecommand{\MRhref}[2]{%
  \href{http://www.ams.org/mathscinet-getitem?mr=#1}{#2}
}
\providecommand{\href}[2]{#2}

\end{document}